\documentclass[12pt, a4paper]{article}

\usepackage[utf8]{inputenc}
\usepackage[french, english]{babel}
\usepackage{indentfirst}
\usepackage[T1]{fontenc}
\usepackage{amsmath, amssymb, amsthm, stmaryrd}
\usepackage{url}
\usepackage{hyperref}
\usepackage{geometry}
\usepackage{color}
\usepackage{dsfont}
\usepackage{todonotes}
\usepackage{enumitem}
\usepackage[mathscr]{euscript}
\usepackage{mathtools}
\usepackage[ddmmyyyy,hhmmss]{datetime}
\usepackage{xcolor}
\usepackage{bbold}

\newcounter{compt}

\newtheorem{theorem}[compt]{Théorème}

\newtheorem{lemme}[compt]{Lemme}

\newcommand\R{\mathbb{R}}
\newcommand\N{\mathbb{N}}

\newcommand\Z{\mathbb{Z}}
\renewcommand\leq{\leqslant}
\renewcommand\geq{\geqslant}
\renewcommand\theta{\vartheta}

\newcommand\e{\mathrm{e}}
\renewcommand\epsilon{\varepsilon}

\begin{document}
\title{Concentration simultanée\\
de fonctions additives}
\author{Élie \bsc{Goudout}}
\date{}
\maketitle
\vspace{-.5cm}
\begin{abstract}
We study the simultaneous concentration of the values of several additive functions along polynomial shifts. Under a slight restriction, this yields an extension of a result from Halász in 1975.
\end{abstract}

\section{Introduction et énoncé des résultats}

Étant donné $f$ une fonction additive et $r\geq 0$ une fonction multiplicative, pour tout $x\geq 2$ on note
\[E_f(x ; r):=1+\sum_{\substack{p\leq x\\ f(p)\neq 0}}\frac{r(p)}{p},\]
et on pose $E_f(x ;1):=E_f(x)$. En 1975, Halász~\cite{halaszconcentration} montre qu'uniformément pour $f$ additive et $x\geq 1$, on a
\[\sup_{k\in\R}\#\left\{n\leq x\,:\quad f(n)=k\right\}\ll\frac{x}{\sqrt{E_f(x)}}.\]
Dans le cas où $f=\omega$ est la fonction nombre de facteurs premiers, lorsque $x$ est grand on a\footnote{Ici et dans la suite, on note $\log_k$ la $k$-ième itérée de la fonction $\log$. ($k\geq 2$)} $E_f(x)=\log_2 x+O(1)$, et la majoration est optimale lorsque l'entier $k$ vérifie $k=\log_2 x+O(\sqrt{\log_2 x})$.

On s'intéresse à une généralisation de ce théorème lorsqu'on fixe plusieurs valeurs de fonctions additives. Il s'agit entre autres, étant données $f$ et $g$ deux fonctions additives, de majorer le cardinal
\[\sup_{k, k'\in\R}\#\left\{n\leq x\,:\quad f(n)=k,\ g(n+1)=k'\right\}.\]
Le cas $f=g=\omega$ a été traité dans~\cite{goudoutpikshort}, puis avec plus de généralité dans~\cite{tenenbaumconcentration}. En effet, y est prouvée l'inégalité
\begin{equation}\label{kjshdfkjshfksueyrzie}
\sup_{k, k'\in\Z}\#\left\{n\leq x\,:\quad \omega(n)=k,\ \omega(n+1)=k'\right\}\ll\frac{x}{\log_2 x}.
\end{equation}
Dans ces deux articles, la majoration est même explicite en $k,k'\ll\log_2 x$. Ici, on traite le cas des fonctions additives ne prenant qu'un nombre fini de valeurs sur les puissances de nombres premiers. On note qu'en considérant
\[\mathcal{A}:=\left\{n\leq x\,:\quad \log_2x-10\sqrt{\log_2x}\leq \omega(n),\omega(n+1)\leq\log_2x+10\sqrt{\log_2x}\right\},\]
qui vérifie $\vert\mathcal{A}\vert\geq x/2$ quand $x$ est assez grand, on obtient avec~(\ref{kjshdfkjshfksueyrzie})
\begin{equation*}
\sup_{k, k'\in\Z}\#\left\{n\leq x\,:\quad \omega(n)=k,\ \omega(n+1)=k'\right\}\asymp\frac{x}{\log_2 x}.
\end{equation*}

On définit maintenant le cadre général d'étude. Étant donné $r\geq 1$ un entier fixé, comme dans~\cite{tenenbaumconcentration}, on considère une famille $\{Q_j\}_{1\leq j\leq r}$ de polynômes irréductibles de $\Z [X]$, deux à deux premiers entre eux et sans diviseur fixe. On pose $Q:=\prod_{1\leq j\leq r}Q_j$. Pour $m\geq 1$, on note $\rho_j(m)$ (\emph{resp.} $\rho_0(m)$) le nombre de racines de $Q_j$ (\emph{resp.} de $Q$) dans $\Z/m\Z$, et $D_j$ (\emph{resp.} $D$) le discriminant de $Q_j$ (\emph{resp.} de $Q$). On note
\begin{align*}
&g_j:=\deg Q_j\quad(0\leq j\leq r),\quad g:=g_0=\sum_{1\leq j\leq r}g_j,\\
&Q(X)=\sum_{0\leq i\leq g}\beta_iX^i,\quad\beta:=\beta_g,\quad\Vert Q\Vert :=\max_{0\leq i\leq g}\vert\beta_i\vert,\\
&\varphi_j(n):=n\prod_{p\vert n}\left(1-\frac{\rho_j(p)}{p}\right)\quad (n\geq 1,\ 0\leq j\leq r),
\end{align*}
où ici et dans la suite, $p$ désigne un nombre premier. On rappelle quelques bornes classiques (voir~\cite{nagell}). Pour $0\leq j\leq r$ et $\nu\geq 1$,
\begin{align}
\rho_j(p^{\nu})&\leq\min(g_jp^{\nu-1}, g_jp^{\nu-\nu/g_j}, p^{\nu-1}\rho_j(p)),&(p\geq 2)\label{all}\\
\rho_j(p^{\nu+1})&\leq\rho_j(p^{\nu})\leq\min(g_j,p-1).&(p\nmid D_j,\ p\geq 2)\label{no}
\end{align}
D'après~\cite[satz~191]{Landau} et \cite[lemme~3.1]{tenenbaumconcentration}, il existe des constantes $M_j, M_j'\geq 1$ et une constante $c=c(g)$ telles que pour $1\leq j\leq r$ et $x\geq 2$,
\begin{align}
&\left\vert\sum_{p\leq x}\frac{\rho_j(p)}{p}-\log_2 x\right\vert\leq M_j\label{mj}\\
&\sum_{p\leq x}\rho_j(p)=\mathrm{li}(x)+O(M_j'+x\exp(-c\sqrt{\log x})).\label{pirho}
\end{align}
On pose $M:=\sum_{1\leq j\leq r}M_j$ et $M':=\sum_{1\leq j\leq r}M_j'$. Étant donnée une fonction additive ou multiplicative définie sur $\N$, on l'étend naturellement à $\Z$ par parité et en fixant arbitrairement sa valeur en $0$ à $0$. 

\begin{theorem}\label{thr}
Soit $r,g,\mathfrak{V}\geq 1$ des entiers et $0<\varepsilon,\delta,\lambda<1$ des réels, tous fixés. Pour $x\geq 1$, on pose $z:=\e^{(\log x)^{1-\lambda}}$. Il existe des constantes $K, c_0>0$ telles qu'uniformément pour $x\geq c_0\Vert Q\Vert^{\delta}$, $x^{\varepsilon}<y\leq x$, et $f_1,\dots ,f_r$ des fonctions additives vérifiant, pour tout $1\leq j\leq r$,
\begin{align}
&\#\left\{f_j(n)\,:\quad n\leq t^u,\ P^-(n)> t\right\}\leq\mathfrak{V}^{u},&(t\geq z,\ u\geq 1)\label{hajzdklqsyugb}\tag{\textasteriskcentered}
\end{align}
on a
\begin{equation}\label{ineg}
\sup_{k_1,\dots, k_r\in\R}\sum_{\substack{x<n\leq x+y \\ f_j(Q_j(n))=k_j\quad (1\leq j\leq r)}}1\ll\left(\frac{\beta D}{\varphi_0(\beta D)}\right)^{K}\frac{\e^{3M+M'/z}y}{\prod_{1\leq j\leq r}\sqrt{E_{f_j}(x ; \rho_j)}},
\end{equation}
La constante implicite, $K$ et $c_0$ dépendent au plus de $r, g, \mathfrak{V}, \varepsilon, \delta$ et $\lambda$.
\end{theorem}
\noindent On note que dans beaucoup de cas d'étude, $M'\ll z$. En particulier, si les polynômes $Q_1,\dots, Q_r$ sont fixés, on obtient
\[\sup_{k_1,\dots, k_r\in\R}\sum_{\substack{x<n\leq x+y \\ f_j(Q_j(n))=k_j\quad (1\leq j\leq r)}}1\ll\frac{y}{\prod_{1\leq j\leq r}\sqrt{E_{f_j}(x ; \rho_j)}}.\]
Par ailleurs, l'hypothèse~(\ref{hajzdklqsyugb}) ne fait pas intervenir les $f_j(p^{\nu})$ pour $p\leq z$, et elle est automatiquement vérifiée si $f_j(p^{\nu})$ ne dépend que de $\nu$ pour $p>z$. Par exemple, cela est vrai pour la fonction nombre de facteurs premiers, avec ou sans multiplicité. On mentionne aussi que
\[\frac{\beta D}{\varphi_0(\beta D)}\ll\left(\frac{\beta D}{\varphi(\beta D)}\right)^g.\]
Il serait intéressant de pouvoir se passer de l'hypothèse~(\ref{hajzdklqsyugb}), qui semble n'être qu'un artefact de la méthode employée.

Pour compléter le résultat, on s'intéresse à une borne inférieure pour le membre de gauche de~(\ref{ineg}).
\begin{theorem}\label{min}
Soit $r,g,\mathfrak{V}\geq 1$ des entiers et $0<\varepsilon,\delta<1$ des réels, tous fixés. On pose $\varepsilon_0:=\varepsilon/(50g)$ et on note~$($\textasteriskcentered$')$ la condition~$($\ref{hajzdklqsyugb}$)$ restreinte aux $t\geq x^{\varepsilon_0}$. Il existe une constante $c_0>0$ telle qu'uniformément pour $x\geq c_0\Vert Q\Vert^{\delta}$, et $f_1,\dots ,f_r$ des fonctions additives vérifiant $($\textasteriskcentered$')$, on a
\begin{equation}\label{dsfjlskjhfd}
\sup_{k_1,\dots k_r\in\R}\sum_{\substack{x<n\leq x+y \\ f_j(Q_j(n))=k_j}}1\gg\frac{\e^{-M}y}{(\log x)^r}\sup_{k_1,\dots k_r\in\R}\prod_{1\leq j\leq r}\sum_{\substack{a_1,\dots,a_r\leq x^{\varepsilon_0} \\ f_j(a_j)=k_j \\ (a_i,a_j)=(a_j, \beta D)=1}}\frac{\rho_j(a_j)\varphi(a_j)}{a_j^2},
\end{equation}
et en particulier, en notant $\omega_y$ la fonction nombre de facteurs premiers distincts inférieurs ou égaux à $y$, uniformément pour $y_1,\dots, y_r\geq 1$,
\begin{equation}\label{vejfhbnqmeflsdj}
\sup_{k_1,\dots k_r\in\N}\sum_{\substack{x<n\leq x+y \\ \omega_{y_j}(Q_j(n))=k_j}}1\gg\left(\frac{\varphi(\beta D)}{\beta D}\right)^g\frac{\e^{-2M}y}{\sqrt{\prod_{1\leq j\leq r}E_{\omega_{y_j}}(x ; \rho_j)}}.
\end{equation}
Les constantes implicites et $c_0$ dépendent au plus de $r, g,\mathfrak{V}, \varepsilon$ et $\delta$
\end{theorem}
Dans le cas où les polynôme $Q_1,\dots,Q_r$ sont fixés et $f_j=\omega_{y_j}$ pour tout $j$, les estimations~(\ref{ineg}) et~(\ref{vejfhbnqmeflsdj}) sont du même ordre de grandeur.

\section{Remerciements}

Je tiens à remercier chaleureusement Régis de la Bretèche pour nos nombreuses discussions, ainsi que ses suggestions toujours éclairées.

\section{Majoration}

Dans cette section, on démontre le Théorème~\ref{thr}. Pour cela, on commence par démontrer deux lemmes.

\begin{lemme}\label{fqhdvgbjhushjk}
Uniformément pour $B_n\geq 0$ $(n\in\Z^*)$ des réels, on a
\[\int_0^1\Big(1+\sum_{n\in\Z^*}\sin^2(\pi nt)B_n\Big)\e^{-\sum_{n\in\Z^*}\sin^2(\pi nt)B_n}\mathrm{d}t\ll\frac{1}{\sqrt{1+\sum_{n\in\Z^*} B_n}}.\]
\end{lemme}

\begin{proof}[Démonstration]
On montre d'abord qu'uniformément pour $B>0$, on a
\[\int_0^1\e^{-\sin^2(\pi t)B}\mathrm{d}t\ll\frac{1}{\sqrt{1+B}}.\]
En posant $\sin(\pi t)\sqrt{B}=u$, on obtient
\[\int_0^1\e^{-\sin^2(\pi t)B}\mathrm{d}t=\frac{2}{\pi}\int_0^{\sqrt{B}}\frac{\e^{-u^2}}{\sqrt{B-u^2}}\mathrm{d}u\ll\frac{1}{\sqrt{1+B}}.\]
Dans le cas général maintenant, quitte à ne sommer que sur les $B_n\neq 0$, on peut supposer que tous les $B_n$ sont non nuls. De même, quitte à ne considérer que les sommes partielles pour $0<\vert n\vert\leq N$ puis à faire tendre $N$ vers $+\infty$, on suppose que $\sum_{n\in\Z^*}B_n<\infty$. On définit alors, pour $n\in\Z^*$, le réel $\theta_n$ tel que $\theta_nB_n=\sum_{i\in\Z*}B_i=:B$. On vérifie que $\sum_{n\in\Z^*}\frac{1}{\theta_n}=1$. Il vient
\begin{align*}
\int_0^1\Big(1+\sum_{n\in\Z^*}\sin^2(\pi nt)&B_n\Big)\e^{-\sum_{n\in\Z^*}\sin^2(\pi nt)B_n}\mathrm{d}t\\
&\leq2\int_0^1\prod_{n\in\Z^*}\Big(\e^{-\frac{1}{2}\sin^2(\pi nt)B_n}\Big)\mathrm{d}t\\
&\leq2\prod_{n\in\Z^*}\Big(\int_0^1\e^{-\sin^2(\pi nt)\frac{\theta_n}{2}B_n}\mathrm{d}t\Big)^{1/\theta_n}\\
&\leq2\prod_{n\in\Z^*}\Big(\int_0^1\e^{-\sin^2(\pi t)\frac{B}{2}}\mathrm{d}t\Big)^{1/\theta_n}\\
&\ll\frac{1}{\sqrt{1+B}},
\end{align*}
où l'on a successivement utilisé $1+x\leq 2\e^{x/2}$ pour $x\geq 0$, l'inégalité de Hölder et la 1-périodicité de $t\mapsto\sin^2(\pi t)$.
\end{proof}

On démontre maintenant une généralisation du théorème de Halász pour une seule fonction additive, dans le cas de poids non constants. Dans la suite, pour $1\leq y\leq x$, on note
\begin{equation}\label{yeghfjbjnsd}
\mathcal{S}(x, y):=\left\{n\leq x\,:\quad P^+(n)\leq y\right\},
\end{equation}
l'ensemble des entiers $y$-friables inférieurs ou égaux à $x$. On utilise aussi la notation classique
\[\mathrm{Li}(t):=\int_2^t\frac{\mathrm{d}t}{\log t}.\hspace{1cm}(t\geq2)\]

\begin{lemme}\label{concentrationuni}
Soient $\varepsilon, b, A>0$ et $0<\lambda<1$ des constantes. Étant donné $y\geq 1$, on pose $z:=\e^{(\log y)^{1-\lambda}}$. Uniformément pour $C,C'\geq 1$, $2\leq y\leq x$ tels que $(\log y)^{\varepsilon(1-\lambda)}\geq (\log_2 x)^2$, $f$ une fonction additive telle que
\begin{equation}\label{condition}
\#\left\{f(p)\,:\quad t<p\leq t^{\e}\right\}\leq A,\hspace{2cm}(z\leq t\leq y)
\end{equation}
et $r\geq 0$ une fonction multiplicative telle que
\begin{gather}
\max_{p\leq y}r(p)\leq A,\hspace{1cm}\sum_{\substack{p\leq y \\ \nu\geq 2}}\frac{r(p^{\nu})\log(p^{\nu})}{p^{\nu}}\leq A,\label{sdkfjh1}\\
\left\vert\sum_{p\leq t}\frac{r(p)}{p}-b\log_2 t\right\vert\leq C,\hspace{1cm}(2\leq t\leq y)\label{sdkfjh2}
\end{gather}
pour laquelle il existe une fonction multiplicative $\widetilde{r}$ telle que
\begin{gather}
\sum_{p\leq y}\frac{\vert r(p)-\widetilde{r}(p)\vert}{p}\leq A,\label{vsdjhfbjsdhbf}\\
\left\vert\sum_{p\leq t}\widetilde{r}(p)-b\mathrm{Li}(t)\right\vert\leq A\left(C'+\frac{t}{(\log t)^{1+\varepsilon}}\right),\hspace{1cm}(2\leq t\leq y)\label{jhbdscjnqhucshwb}
\end{gather}
on a
\[\sup_{k\in\R}\sum_{\substack{n\in\mathcal{S}(x, y) \\ f(n)=k}}r(n)\ll\frac{(\log y)^b}{\log x}\frac{\e^{2C+C'/z}x}{\sqrt{E_{f}(y ; r)}},\]
où la constante implicite dépend au plus de $\varepsilon, b, A$ et $\lambda$.
\end{lemme}
\noindent On note que la condition~(\ref{condition}) peut être affaiblie en
\[\#\left\{f(p)\,:\quad t<p\leq t^{1+c}\right\}\leq A\hspace{2cm}(z\leq t\leq y)\]
avec $c:=(\log_2 x)^2/(\log y)^{\varepsilon(1-\lambda)}\leq 1$, par adaptation directe de la méthode. Il est par ailleurs possible de démontrer une version uniforme en $\varepsilon$.

\begin{proof}[Démonstration]
Comme le fait remarquer Halász au début de la démonstration du théorème de~\cite{halaszconcentration}, quitte à modifier $f$ d'une manière précise, on peut supposer qu'elle est à valeurs entières. La construction qu'il emploie garantit que le nombre de valeurs prises par $f(p)$ n'augmente pas, et que $E_{f}(y ; r)$ ne diminue pas. En posant $R(t):=\sum_{\substack{n\in\mathcal{S}(x, y)}}r(n)\e^{2i\pi f(n)t}$, pour tout $k\in\Z$, on a alors
\begin{equation*}
\sum_{\substack{n\in\mathcal{S}(x, y) \\ f(n)=k}}r(n)=\int_0^{1}R(t)\e^{-2ik\pi t}\mathrm{d}t\leq\int_0^{1}\vert R(t)\vert\mathrm{d}t.
\end{equation*}
En utilisant~(\ref{sdkfjh1}-\ref{sdkfjh2}), on applique le théorème~1.1 de~\cite{tenenbaummoyenneseffectives} avec $T:=(\log_2 x)^2$ et on utilise la démonstration du corollaire~2.1 de ce même article, pour obtenir uniformément pour $t\in[0, 1]$,
\begin{equation}\label{fsghdjkfqjs}
\vert R(t)\vert\ll\e^{C}x\frac{(\log y)^b}{\log x}\left\{\frac{1+m(t)}{\e^{m(t)}}+\frac{1}{\log_2 x}\right\}
\end{equation}
où
\[m(t):=\min_{\vert\tau\vert\leq (\log_2 x)^2}\sum_{p\leq y}\frac{r(p)(1-\cos(2\pi f(p)t-\tau\log p))}{p}.\]
On pose, pour $k\geq 0$ et $t\in[0, 1]$,
\begin{align*}
\gamma_{k,t}(u)&:=1-\max_{\substack{k<\log_2 p\leq k+1 \\ z<p\leq y}}\cos(2\pi f(p)t-u),\hspace{1cm}(u\in\R)\\
s_t&:=\min_{k\geq 0}\frac{1}{2\pi}\int_0^{2\pi}\gamma_{k,t}(u)\mathrm{d}u\geq 1-\frac{A}{\pi}\sin\left(\frac{\pi}{A}\right)\gg 1,
\end{align*}
avec la convention $\max_{\emptyset}=0$. Ainsi, pour tous $z\leq y_0\leq y$, on a pour un certain $\vert\tau\vert\leq(\log_2 x)^2$,
\begin{align*}
m(t)&\geq\sum_{k\geq 0}\sum_{\substack{k<\log_2 p\leq k+1 \\ y_0<p\leq y}}\frac{r(p)\gamma_{k,t}(\tau\log p)}{p}\\
&\geq\sum_{k\geq 0}\sum_{\substack{k<\log_2 p\leq k+1 \\ y_0<p\leq y}}\frac{\widetilde{r}(p)\gamma_{k,t}(\tau\log p)}{p}-A,
\end{align*}
où l'on a utilisé~(\ref{vsdjhfbjsdhbf}).
De manière analogue au lemme 4.13 de~\cite[III.4]{Tenenbaum}, on estime la somme ci-dessus avec~(\ref{jhbdscjnqhucshwb}), par sommation par parties. On obtient ainsi, uniformément pour $z\leq y_0\leq y$ et $t\in[0,1]$, lorsque $\tau\neq 0$,
\begin{align*}
m(t)&\geq s_tb\log\left(\frac{\log y}{\log y_0}\right)\\
+&O\Bigg(\sum_{\log_2 y_0-1<k\leq\log_2 y}\frac{1}{\vert\tau\vert\e^k}+(1+\vert\tau\vert)\left(\frac{1}{(\e^k)^{\varepsilon}}+\frac{C'}{\e^{\e^k}}\right)\Bigg)\\
&\geq s_tb\log\left(\frac{\log y}{\log y_0}\right)+O\left(\frac{1}{\vert\tau\vert\log y_0}+(1+\vert\tau\vert)\left(\frac{1}{(\log y_0)^{\varepsilon}}+\frac{C'}{y_0^{1/\e}}\right)\right).
\end{align*}
Si $1\leq\vert\tau\vert\leq(\log_2 x)^2$, on pose $y_0:=z^{3\e}$. Puisque $(\log y)^{\varepsilon(1-\lambda)}\geq(\log_2 x)^2$, on obtient
\[m(t)\geq \lambda s_tb\log_2 y+O\left(1+\frac{C'}{z^2}\right).\]
On définit ensuite $w$ tel que
\[\log w=2\e(\log y)\exp\Bigg(-\frac{\lambda}{b}\bigg\{\sum_{v\in\Z^*}\sin^2(\pi vt)\sum_{\substack{p\leq y \\ f(p)=v}}\frac{r(p)}{p}-2C\bigg\}\Bigg).\]
D'après~(\ref{sdkfjh2}), on a $w\geq z^{2\e}$. Si $w\geq y$, on retient trivialement que $m(t)\geq 0$. Dans le cas contraire, si $\frac{1}{\log w}<\vert\tau\vert\leq 1$, avec $y_0:=w$, on obtient
\[m(t)\geq \lambda s_t\bigg\{\sum_{v\in\Z^*}\sin^2(\pi vt)\sum_{\substack{p\leq y \\ f(p)=v}}\frac{r(p)}{p}-2C\bigg\}+O\left(1+\frac{C'}{z^2}\right).\]
Enfin, lorsque $\vert\tau\vert\leq\frac{1}{\log w}$, on minore trivialement
\begin{align*}
m(t)&\geq\sum_{p\leq w}\frac{r(p)(1-\cos(2\pi f(p)t))}{p}+O(1)\\
&=2\sum_{v\in\Z^*}\sin^2(\pi vt)\sum_{\substack{p\leq w \\ f(p)=v}}\frac{r(p)}{p}+O(1)\\
&\geq 2\sum_{v\in\Z^*}\sin^2(\pi vt)\sum_{\substack{p\leq y \\ f(p)=v}}\frac{r(p)}{p}-2b\log\left(\frac{\log y}{\log w}\right)-4C+O(1)\\
&\geq 2(1-\lambda)\bigg\{\sum_{v\in\Z^*}\sin^2(\pi vt)\sum_{\substack{p\leq y \\ f(p)=v}}\frac{r(p)}{p}-2C\bigg\}+O(1).
\end{align*}
Ainsi, quelle que soit la valeur de $\tau$, en posant $\eta:=\min\left(\lambda s_tb, \lambda s_t, 2(1-\lambda), 1/2\right)$, lorsque $x$ (et donc $z$) est suffisamment grand, on a
\begin{equation*}
m(t)\geq\eta\sum_{v\in\Z^*}\sin^2(\pi vt)\sum_{\substack{p\leq y \\ f(p)=v}}\frac{r(p)}{p}-\Big(C+\frac{C'}{z}+O(1)\Big)
\end{equation*}
et $m(t)\geq 0$. Ainsi, avec~(\ref{fsghdjkfqjs}), on obtient le résultat désiré en intégrant selon~$t$, par application directe du Lemme~\ref{fqhdvgbjhushjk}.
\end{proof}

On démontre désormais le théorème. Soit $2\leq y\leq x$ et $n\in(x,2x]$. Un élément clé de la preuve de~(\ref{kjshdfkjshfksueyrzie}) est que si l'on note $n_y$ la partie $y$-friable de $n$, alors $\omega(n)-\omega(n_y)$ ne peut prendre au plus qu'un nombre fini de valeurs lorsque $\log y\asymp\log x$. Cela est toujours vrai lorsqu'on remplace $\omega$ par une fonction additive vérifiant~(\ref{hajzdklqsyugb}). C'est la seule hypothèse du théorème sur les fonctions additives, et il serait intéressant de pouvoir s'en passer. La démonstration diffère peu de celle de~\cite{tenenbaumconcentration}.

\begin{proof}[Démonstration du Théorème~\ref{thr}]
On suppose donnés les paramètres de l'énoncé. Durant la démonstration, on note $K$ une constante positive qui sera toujours prise suffisamment grande au besoin. Puisque $Q$ n'admet qu'un nombre fini de racines dans $\Z$, on  peut, sans perte de généralité, faire tous nos raisonnements en supposant $Q(n)\neq0$. On se donne $k_1, \dots, k_r\in\R$ fixés. Pour tout entier $n\geq 1$ avec $Q(n)\neq 0$, on note $\xi_n$ le plus grand entier tel que la partie $\xi_n$-friable de $Q(n)$ soit inférieure ou égale à $x^{2\varepsilon/3}$. On obtient ainsi la décomposition canonique $Q(n)=b_n\prod_{1\leq j\leq r}a_{j, n}$ où $b_n$ est la partie $\xi_n$-criblée de $Q(n)$, et pour tout $1\leq j\leq r$, $a_{j, n}\vert Q_j(n)$. On note $p_n:=P^-(b_n)$ et $\nu_n$ tel que $p_n^{\nu_n}\Vert Q(n)$.

On note $N_1(x)$ le nombre d'entiers $n$ apparaissant dans la somme de~(\ref{ineg}) pour lesquels $a_{1, n}\cdots a_{r, n}\leq x^{\varepsilon/3}$. On note $N_2(x)$ le nombre des autres entiers $n$ de la somme de~(\ref{ineg}). Pour la suite, on pose
\[\gamma:=g+1/\delta+1,\]
de sorte que $\vert Q(n)\vert\leq x^{\gamma}$ pour $x$ suffisamment grand et $x<n\leq x+y$ ; et on note
\begin{align*}
&\mathcal{V}_u:=\left\{f_j(n)\,:\quad 1\leq j\leq r,\ n\leq x^{\gamma},\ P^-(n)>x^{\gamma/u}\right\}.&(u\geq 1)
\end{align*}
D'après~(\ref{hajzdklqsyugb}), on a $\vert\mathcal{V}_u\vert\leq r\mathfrak{V}^u$ dès que $x^{\gamma/u}\geq z$.

Considérons $n$ compté dans $N_1(x)$. Alors $p_n^{\nu_n}>x^{\varepsilon/3}$. Par ailleurs, si $p_n>x^{\varepsilon/(6g)}$, en posant $\eta_{_0}:=6g\gamma/\varepsilon$, on a $f_j(Q_j(n))-f_j(a_{j, n})\in\mathcal{V}_{\eta_{_0}}$ pour tout $1\leq j\leq r$ lorsque $x$ est assez grand. Ainsi,
\begin{equation}\label{deux}
N_1(x)\leq\sum_{\substack{a_1\cdots a_r\leq x^{\varepsilon/3} \\ k_j-f_j(a_j)\in \mathcal{V}_{\eta_{_0}}}}\sum_{\substack{x<n\leq x+y \\ a_j\vert Q_j(n) \\ P^-(\frac{Q(n)}{a_1\cdots a_r})>x^{\varepsilon/(6g)}}}1+\sum_{\substack{p\leq x^{\varepsilon/(6g)} \\ p^{\nu}>x^{\varepsilon/3}}}\sum_{\substack{x<n\leq x+y \\ p^{\nu}\Vert Q(n)}}1.
\end{equation}
On majore en premier la seconde double somme. Pour les $p$ considérés, quitte à remplacer $\nu$ par $\nu'\leq\nu$ le plus petit tel que $p^{\nu'}>x^{\varepsilon/3}$ et $p^{\nu}\Vert Q(n)$ par $p^{\nu'}\vert Q(n)$ --- auquel cas $p^{\nu'-1}\leq x^{\varepsilon/3}$ et donc $p^{\nu'}\leq px^{\varepsilon/3}\leq x^{\varepsilon/2}\leq y$ ---, on peut supposer $p^{\nu}\leq y$. On obtient
\[\sum_{\substack{p\leq x^{\varepsilon/(6g)} \\ p^{\nu}>x^{\varepsilon/3}}}\sum_{\substack{x<n\leq x+y \\ p^{\nu}\Vert Q(n)}}1
\leq\sum_{\substack{p\leq x^{\varepsilon/(6g)} \\ x^{\varepsilon/3}<p^{\nu}\leq y}}\sum_{\substack{x<n\leq x+y \\ p^{\nu}\vert Q(n)}}1
\leq 2\sum_{\substack{p\leq x^{\varepsilon/(6g)} \\ p^{\nu}>x^{\varepsilon/3}}}\frac{y\rho_0(p^{\nu})}{p^{\nu}}\ll\frac{y}{x^{\varepsilon/(6g)}},\]
où la dernière inégalité est une conséquence directe de~(\ref{all}).

On majore maintenant la première double somme de~(\ref{deux}). Comme dans~\cite{tenenbaumconcentration}, pour tout $1\leq j\leq r$, on décompose $a_j$ sous la forme $a_j=t_jd_j$ où $t_j\vert(\beta D)^{\infty}$ et $(d_j, \beta D)=1$, de sorte que $k_j-f_j(t_j)-f_j(d_j)\in\mathcal{V}_{\eta_{_0}}$ et $(d_i, d_j)=1$ pour tous $i\neq j$. En posant $T:=t_1\cdots t_r$, on a alors
\begin{align*}
\sum_{\substack{a_1\cdots a_r\leq x^{\varepsilon/3} \\ k_j-f_j(a_j)\in\mathcal{V}_{\eta_{_0}}}}\sum_{\substack{x<n\leq x+y \\ a_j\vert Q_j(n) \\ P^-(\frac{Q(n)}{a_1\cdots a_r})>x^{\varepsilon/(6g)}}}1
&\ll\sum_{\substack{t_1d_1\cdots t_rd_r\leq x^{\varepsilon/3} \\ T\vert (\beta D)^{\infty} \\ (d_i,d_j)=(d_j,\beta D)=1 \\ k_j-f_j(t_j)-f_j(d_j)\in\mathcal{V}_{\eta_{_0}}}}\sum_{\substack{x<n\leq x+y \\ T\vert Q(n) \\ d_j\vert Q_j(n) \\ p\vert Q(n)\Rightarrow p\vert Td_1\cdots d_r\text{ ou }p>x^{\varepsilon/(6g)}}}1\\
\end{align*}
Avec le même raisonnement de crible que~\cite{tenenbaumconcentration}, on majore la dernière somme, sous les hypothèses $Td_1\cdots d_r\leq x^{\varepsilon/3}$, $T\vert(\beta D)^{\infty}$ et $(d_i, d_j)=(d_j,\beta D)=1$ :
\[\sum_{\substack{x<n\leq x+y \\ T\vert Q(n) \\ d_j\vert Q_j(n) \\ p\vert Q(n)\Rightarrow p\vert Td_1\cdots d_r\text{ ou }p>x^{\varepsilon/(6g)}}}1\ll \frac{\e^My}{(\log x)^r}\left(\frac{\beta D}{\varphi_0(\beta D)}\right)^r\frac{\rho_0(T)}{T}\prod_{1\leq j\leq r}\frac{\rho_j(d_j)}{\varphi_j(d_j)}.\]
En posant
\begin{align*}
H&:=\e^M\left(\frac{\beta D}{\varphi_0(\beta D)}\right)^r,\\
S_j(t_j)&:=\sum_{\substack{d_j\leq x \\ (d_j,\beta D)=1 \\ k_j-f_j(t_j)-f_j(d_j)\in\mathcal{V}_{\eta_{_0}}}}\frac{\rho_j(d_j)}{\varphi_j(d_j)},&(1\leq j\leq r)
\end{align*}
on obtient alors
\[N_1(x)\ll\frac{Hy}{(\log x)^r}\sum_{\substack{t_1\cdots t_r\vert (\beta D)^{\infty}}}\frac{\rho_0(T)}{T}\prod_{1\leq j\leq r}S_j(t_j)+\frac{y}{x^{\varepsilon/(6g)}}.\]

Lorsque $x$ est assez grand, on a $\vert\mathcal{V}_{\eta_{_0}}\vert\leq r\mathfrak{V}^{\eta_{_0}}\ll 1$. Pour $1\leq j\leq r$, on définit la fonction multiplicative $\tilde{\rho}_j$ par
\begin{align*}
\tilde{\rho}_j(p^{\nu})&:=\frac{\rho_j(p^{\nu})p^{\nu}}{\varphi_j(p^{\nu})},&(p\nmid\beta D,\ \nu\geq 1)\\
\tilde{\rho}_j(p^{\nu})&:=0,&(p\vert\beta D,\ \nu\geq 2)\\
\tilde{\rho}_j(p)&:=\rho_j(p),&(p\vert\beta D)
\end{align*}
de sorte que
\[S_j(t_j)\ll\sup_{k\in\R}\sum_{\substack{n\leq x \\ f_j(n)=k}}\frac{\tilde{\rho}_j(n)}{n}.\]
Puisque $\rho_j(p)\leq\min(p-1, g_j)$, on a 
\begin{align}
&\sum_{p\leq t}\frac{\tilde{\rho}_j(p)}{p}=\sum_{p\leq t}\frac{\rho_j(p)}{p}+O(1).&(t\geq 2)\label{azjfhkzbjsn}
\end{align}
On pose
\begin{equation}\label{thfjnlzesjn}
x_0:=\e^{(\log x)^{1-\lambda/2}}.
\end{equation}
Soit $\lambda'>0$ tel que $(1-\lambda/2)(1-\lambda')\geq 1-\lambda$. Avec~(\ref{mj}-\ref{pirho}) et~(\ref{azjfhkzbjsn}), on applique le Lemme~\ref{concentrationuni} à $\tilde{\rho}_j$. Ainsi, pour tout $x'\in(x_0, x]$, en posant $z':=\e^{(\log x')^{1-\lambda'}}$, on obtient
\[\sup_{k\in\R}\sum_{\substack{n\leq x' \\ f_j(n)=k}}\tilde{\rho}_j(n)\ll\frac{\e^{2M_j+M_j'/z'}x'}{\sqrt{E_{f_j}(x' ; \rho_j)}},\]
puisque $E_{f_j}(x' ;\tilde{\rho}_j)\asymp E_{f_j}(x' ;\rho_j)$. Par ailleurs, on a facilement
\[\sum_{n\leq x_0}\frac{\tilde{\rho}_j(n)}{n}\ll\e^{M_j}\log x_0=\e^{M_j}(\log x)^{1-\lambda/2}.\]
Finalement, via une sommation par parties sur $(x_0, x]$, on obtient pour tout $1\leq j\leq r$,
\begin{equation}
S_j(t_j)\ll\frac{\e^{2M_j+M_j'/z}\log x}{\sqrt{E_{f_j}(x ; \rho_j)}}.\label{tghjaifjqisgcbj}
\end{equation}
La majoration voulue pour $N_1(x)$ découle alors de l'inégalité
\begin{align*}
\sum_{\substack{T\vert (\beta D)^{\infty}}}\frac{\rho_0(T)\tau_r(T)}{T}\ll\left(\frac{\beta D}{\varphi_0(\beta D)}\right)^{K},
\end{align*}
qui est une conséquence facile de~(\ref{all}).

On majore désormais $N_2(x)$. Pour cela, on introduit $q_n:=P^+(a_{1, n}\cdots a_{r, n})$ et on pose $\eta(q_n):=\gamma(\log x)/\log q_n$, de sorte que $f_j(Q_j(n))-f_j(a_{j, n})\in\mathcal{V}_{\eta(q_n)}$ pour tout $1\leq j\leq r$ lorsque $x$ est assez grand. Avec~(\ref{hajzdklqsyugb}), pour $q_n\geq z$, on a $\vert\mathcal{V}_{\eta(q_n)}\vert\leq r\mathfrak{V}^{\eta(q_n)}$. En décomposant les $a_j$ comme dans le cas de $N_1(x)$, on obtient
\begin{align}
N_2(x)&\leq\sum_{q\leq x^{2\varepsilon/3}}\sum_{\substack{x^{\varepsilon/3}<t_1d_1\cdots t_rd_r\leq x^{2\varepsilon/3} \\ P^+(Td_1\cdots d_r)=q \\ T\vert(\beta D)^{\infty} \\ (d_i, d_j)=(d_j,\beta D)=1 \\ k_j-f_j(t_j)-f_j(d_j)\in\mathcal{V}_{\eta(q)}}}\sum_{\substack{x<n\leq x+y \\ T\vert Q(n) \\ d_j\vert Q_j(n) \\ p\vert Q(n)\Rightarrow p\vert Td_1\cdots d_r\text{ ou }p>q}}1.\label{jhgsjfdh}
\end{align}
Comme précédemment, sous les hypothèses $Td_1\cdots d_r\leq x^{2\varepsilon/3}$, $T\vert(\beta D)^{\infty}$ et $(d_i, d_j)=(d_j,\beta D)=1$, on majore la dernière somme
\[\sum_{\substack{x<n\leq x+y \\ T\vert Q(n) \\ d_j\vert Q_j(n) \\ p\vert Q(n)\Rightarrow p\vert Td_1\cdots d_r\text{ ou }p>q}}1\ll\frac{Hy}{(\log q)^r}\frac{\rho_0(T)}{T}\prod_{1\leq j\leq r}\frac{\rho_j(d_j)}{\varphi_j(d_j)}.\]
Soit $C>0$ une constante que l'on fixera plus tard. Pour $q\geq 2$, on pose
\[\alpha=\alpha(q):=C/\log q.\]
Dans~(\ref{jhgsjfdh}), on majore trivialement la contribution des $q\leq\e^{2gC}$. Pour cela, on remarque qu'avec~(\ref{all}), lorsque $P^+(m)\leq\e^{2gC}$ et $0\leq j\leq r$, on a $\varphi_j(m)\gg m$ et $\rho_j(m)\ll m^{1-1/g_j}$. Pour le reste, on utilise l'astuce de Rankin en introduisant $(Td_1\cdots d_r/x^{\varepsilon/3})^{\alpha}>1$. On obtient alors
\begin{equation}\label{fghjkazdu}
N_2(x)\ll Hy\sum_{\e^{2gC}<q\leq x^{2\varepsilon/3}}\frac{1}{(\log q)^rx^{\varepsilon\alpha/3}}\hspace{-.9cm}\sum_{\substack{t_1d_1\cdots t_rd_r\leq x^{2\varepsilon/3} \\ P^+(Td_1\cdots d_r)=q \\ T\vert(\beta D)^{\infty} \\ (d_j,\beta D)=1 \\ k_j-f_j(t_j)f_j(d_j)\in\mathcal{V}_{\eta(q)}}}\hspace{-.9cm}\frac{\rho_0(T)}{T^{1-\alpha}}\prod_{1\leq j\leq r}\frac{\rho_j(d_j)d_j^{\alpha}}{\varphi_j(d_j)}+\frac{Hy}{x^{\varepsilon/(4g)}}.
\end{equation}
Pour $\e^{2gC}< q\leq x_0$, on a $\alpha\leq 1/(2g)$ et avec~(\ref{thfjnlzesjn}), $x^{\varepsilon\alpha/3}\gg\e^{2(\log x)^{\lambda/3}}$. On majore alors la contribution des $q\in(\e^{2gC}, x_0]$ avec~(\ref{all}-\ref{no}) en négligeant la condition sur les valeurs de $f_j(d_j)$. Pour cela, on remarque notamment qu'uniformément pour $p\leq q$ et $\nu\geq 1$,
\begin{align}
&\frac{\rho_0(p^{\nu})}{p^{\nu(1-\alpha)}}\ll\min\left(\frac{1}{p}, \frac{1}{p^{\nu/(2g)}}\right),\label{ghjfhajeds}\\
&\frac{\rho_j(p^{\nu})p^{\nu\alpha}}{\varphi_j(p^{\nu})}\ll\frac{1}{p^{\nu(1-\alpha)}}.&(p\nmid\beta D,\ 1\leq j\leq r)\label{thezfdslquhk}
\end{align}
On obtient alors
\begin{multline*}
\sum_{\substack{\e^{2gC}<q\leq x_0}}\frac{1}{(\log q)^rx^{\varepsilon\alpha/3}}\sum_{\substack{Td_1\cdots d_r\leq x^{2\varepsilon/3} \\ P^+(Td_1\cdots d_r)=q \\ T\vert(\beta D)^{\infty} \\ (d_j,\beta D)=1}}\frac{\rho_0(T)\tau_r(T)}{T^{1-\alpha}}\prod_{1\leq j\leq r}\frac{\rho_j(d_j)d_j^{\alpha}}{\varphi_j(d_j)}\\
\ll \left(\frac{\beta D}{\varphi_0(\beta D)}\right)^{K}\e^{-(\log x)^{\lambda/3}},
\end{multline*}
Finalement, il suffit de montrer que pour $q\in(x_0, x^{2\varepsilon/3}]$, la somme intérieure de~(\ref{fghjkazdu}) est
\[\ll\left(\frac{\beta D}{\varphi_0(\beta D)}\right)^{K}\frac{\e^{2M+M'/z}(\log x)^r\mathfrak{V}^{r\eta(q)}}{q\prod_{1\leq j\leq r}\sqrt{E_{f_j}(q ;\rho_j)}}.\]
Pour cela, on montre que pour $1\leq j\leq r$,
\begin{align*}
&\sup_{k\in\R}\sum_{\substack{d\in\mathcal{S}(x, q) \\ (d,\beta D)=1 \\ f_j(d)=k}}\frac{\rho_j(d)d^{\alpha}}{\varphi_j(d)}\ll\frac{\e^{2M+M'/z}(\log x)}{\sqrt{E_{f_j}(q ;\rho_j)}},\\
\sum_{\nu\geq 1}&\sup_{k\in\R}\sum_{\substack{d\in\mathcal{S}(x, q-1) \\ (q^{\nu}d,\beta D)=1 \\ f_j(q^{\nu}d)=k}}\frac{\rho_j(q^{\nu}d)(q^{\nu}d)^{\alpha}}{\varphi_j(q^{\nu}d)}\ll\frac{\e^{2M+M'/z}(\log x)}{q\sqrt{E_{f_j}(q ;\rho_j)}}.
\end{align*}
Avec~(\ref{thezfdslquhk}), la première inégalité implique facilement la seconde. De manière analogue à~(\ref{tghjaifjqisgcbj}), la première inégalité est conséquence du Lemme~\ref{concentrationuni}. Afin d'en vérifier les hypothèses, on remarque que
\[\sum_{p\leq q}\frac{\rho_j(p)}{p^{1-\alpha}}=\sum_{p\leq q}\frac{\rho_j(p)}{p}+O(1).\]
Par ailleurs, de manière analogue au cas $\e^{2gC}<q\leq x_0$, la somme sur $T$ se majore avec~(\ref{ghjfhajeds}). Pour conclure, il suffit de noter que l'on a
\[\sum_{x_0< q\leq x^{2\varepsilon/3}}\frac{(\log x)^r\mathfrak{V}^{r\eta(q)}}{q(\log q)^rx^{\varepsilon\alpha/3}\prod_{1\leq j\leq r}\sqrt{E_{f_j}(q ;\rho_j)}}\ll\frac{1}{\prod_{1\leq j\leq r}\sqrt{E_{f_j}(x ; \rho_j)}}\]
dès que $\varepsilon C>3r\gamma\log(\mathfrak{V})$.
\end{proof}

\section{Minoration}

\begin{proof}[Démonstration du Théorème~\ref{min}]
On suppose donnés les paramètres de l'énoncé. Pour tout $x<n\leq x+y$ tel que $Q(n)\neq 0$ et $1\leq j\leq r$, on pose $a_{j, n}$ la partie $x^{\varepsilon_0}$-friable de $Q_j(n)$. Pour $x$ suffisamment grand, on a $\vert Q(n)\vert\leq x^{g+1/\delta+1}$. Ainsi, puisque $f_1,\dots f_r$ vérifient~(\ref{hajzdklqsyugb}) pour $t\geq x^{\varepsilon_0}$, $f_j(Q_j(n))-f_j(a_{j, n})$ ne peut prendre qu'un nombre fini de valeurs lorsque $x<n\leq x+y$. On peut donc, sans perte de généralité, supposer que pour tout $p>x^{\varepsilon_0}$, on a $f_j(p^{\nu})=0$ lorsque $\nu\geq 1$ et $1\leq j\leq r$. Dans ce cas-là, pour tout $k_1,\dots k_r\in\R$, on a
\[\sum_{\substack{x<n\leq x+y \\ f_j(Q_j(n))=k_j}}1\geq\sum_{\substack{a_1\cdots a_r\in\mathcal{S}(x,x^{\varepsilon_0}) \\ f_j(a_j)=k_j}}\sum_{\substack{x<n\leq x+y \\ a_j\vert Q_j(n) \\ P^-\left(\frac{Q(n)}{a_1\cdots a_r}\right)>x^{\varepsilon_0}}}1,\]
où l'on a utilisé la notation~(\ref{yeghfjbjnsd}).

On introduit quelques notations utilisées dans~\cite{henriot} et~\cite{henrioterr}. Pour $n\geq 1$, on note $\kappa(n)$ le noyau sans facteurs carrés de $n$. Si de plus $m\geq 1$ est un autre entier, on écrit $n\Vert m$ pour signifier $n\vert m$ et $(n, \frac{m}{n})=1$. Par ailleurs, pour $a_1,\dots, a_r\geq 1$, on pose
\begin{multline*}
\check{\rho}(a_1,\dots,a_r):=\#\left\{n\mod{[a_1\kappa(a_1)\cdots a_r\kappa(a_r)]}\,:\right.\\
\left.\forall 1\leq j\leq r,\ a_j\Vert Q_j(n),\ a_1\cdots a_r\Vert Q(n)\right\}.
\end{multline*}
D'après la version corrigée dans~\cite{henrioterr} de~\cite[lemma~6]{henriot}, on a alors
\[\sum_{\substack{x<n\leq x+y \\ f_j(Q_j(n))=k_j}}1\gg y\sum_{\substack{a_1\cdots a_r\in\mathcal{S}(x^{3\varepsilon/25},x^{\varepsilon_0}) \\ f_j(a_j)=k_j}}\frac{\check{\rho}(a_1,\dots,a_r)}{[a_1\kappa(a_1)\cdots a_r\kappa(a_r)]}\prod_{\substack{g<p\leq x^{\varepsilon_0} \\ p\nmid a_1\cdots a_r}}\left(1-\frac{\rho_0(p)}{p}\right).\]
Puisque $\rho_0(p)\leq \sum_{1\leq j\leq r}\rho_j(p)$, avec~(\ref{all}) et~(\ref{mj}), on a
\begin{equation*}
\prod_{\substack{g<p\leq x^{\varepsilon_0} \\ p\nmid a_1\cdots a_r}}\left(1-\frac{\rho_0(p)}{p}\right)\gg\frac{\e^{-M}}{(\log x)^r}.
\end{equation*}
Pour obtenir~(\ref{dsfjlskjhfd}), il suffit de montrer que si l'on impose $(a_i,a_j)=(a_j,\beta D)=1$ pour tous $1\leq i<j\leq r$, alors
\begin{equation}\label{hgefjhqsjnlc}
\frac{\check{\rho}(a_1,\dots,a_r)}{[a_1\kappa(a_1)\cdots a_r\kappa(a_r)]}\geq\prod_{1\leq j\leq r}\frac{\rho_j(a_j)\varphi(a_j)}{a_j^2}.
\end{equation}
Pour cela, on commence par noter comme dans la démonstration du Théorème~\ref{thr}, que lorsque $p\nmid \beta D$, alors pour tous $i\neq j$, on a $p\nmid (Q_i(n), Q_j(n))$. Donc lorsque $p\nmid\beta D$, si $p^{\nu}\Vert Q_j(n)$ pour un certain $j$, alors $p^{\nu}\Vert Q(n)$. La minoration~(\ref{hgefjhqsjnlc}) découle alors de la multiplicativité de $\check{\rho}$ et de~(\ref{no}).

Il reste à démontrer~(\ref{vejfhbnqmeflsdj}). Pour cela, on utilise entre autres le \og $W$-trick\fg\ afin de gérer les conditions $(a_i, a_j)=1$. Soit donc $w\geq 1$ un entier, que l'on fixera plus tard, et $W:=\prod_{p\leq w}p$. On se donne par ailleurs $C\geq 1$, une autre quantité à fixer plus tard. Soit $y_1,\dots, y_r\geq 2$. Pour $1\leq j\leq r$, on introduit $y_j^*:=\min(y_j, x^{\varepsilon_0/C})$ et on pose
\begin{align*}
L_j(\alpha)&:=\sum_{\substack{w<p\leq y_j^* \\ p\nmid\beta D}}\frac{\rho_j(p)\varphi(p)p^{\alpha}}{p^2}, &(0\leq\alpha\leq(\log x)^{-\varepsilon_0/C})\\
L_j&:=\sum_{\substack{w<p\leq y_j^* \\ p\nmid\beta D}}\frac{\rho_j(p)}{p}.
\end{align*}
Puisque $0\leq\rho_j(p)\leq g_j$ pour tout $1\leq j\leq r$, une sommation par parties fournit
\begin{equation}\label{thzfieosdubk}
L_j(\alpha)=L_j(0)+O(1)=L_j+O(1).
\end{equation}
On démontre que la minoration~(\ref{vejfhbnqmeflsdj}) est valable pour les valeurs particulières
\begin{equation}
k_j:=\lfloor L_j\rfloor.\hspace{2cm}(1\leq j\leq r)\label{jhljshgdflh}
\end{equation}
Pour cela, on commence par démontrer le lemme suivant.
\begin{lemme}\label{hejfsndlsdfuh}
Avec les notations ci-dessus, lorsque $w$ et $C$ sont suffisamment grands, pour tout $1\leq j\leq r$, uniformément pour $0\leq k_j'\leq k_j$ on a
\[\sum_{\substack{n\leq x^{\varepsilon_0} \\ \omega_{y_j}(n)=k_j' \\ (n, W\beta D)=1 \\ \mu(n)^2=1}}\frac{\rho_j(n)\varphi(n)}{n^2}\asymp\frac{L_j^{k_j'}}{k_j'!}\exp\bigg\{\sum_{\substack{y_j^*<p\leq x \\ p\nmid\beta D}}\frac{\rho_j(p)}{p}\bigg\}.\]
\end{lemme}
\begin{proof}[Démonstration du lemme]
La majoration est aisée en remplaçant la condition $n\leq x^{\varepsilon_0}$ par $P^+(n)\leq x$. En effet, on obtient
\begin{align*}
\sum_{\substack{n\leq x^{\varepsilon_0} \\ \omega_{y_j}(n)=k_j' \\ (n, W\beta D)=1 \\ \mu(n)^2=1}}\frac{\rho_j(n)\varphi(n)}{n^2}&\leq\frac{L_j(0)^{k_j'}}{k_j'!}\prod_{\substack{y_j^*<p\leq x \\ p\nmid\beta D}}\bigg\{1+\frac{\rho_j(p)}{p}\bigg\}\\
&\ll\frac{L_j^{k_j'}}{k_j'!}\exp\bigg\{\sum_{\substack{y_j^*<p\leq x \\ p\nmid\beta D}}\frac{\rho_j(p)}{p}\bigg\},
\end{align*}
où l'on a utilisé~(\ref{thzfieosdubk}). On démontre maintenant la borne inférieure. Pour cela, on introduit les notations
\[\Sigma_1:=\sum_{\substack{P^+(n)\leq x^{\varepsilon_0/C} \\ \omega_{y_j}(n)=k_j' \\ (n,W\beta D)=1 \\ \mu(n)^2=1}}\frac{\rho_j(n)\varphi(n)}{n^2},\hspace{1cm}\Sigma_2:=\sum_{\substack{P^+(n)\leq x^{\varepsilon_0/C} \\ n>x^{\varepsilon_0} \\ \omega_{y_j}(n)=k_j' \\ (n,W\beta D)=1 \\ \mu(n)^2=1}}\frac{\rho_j(n)\varphi(n)}{n^2},\]
de sorte que l'on a
\[\sum_{\substack{n\leq x^{\varepsilon_0} \\ \omega_{y_j}(n)=k_j' \\ (n, W\beta D)=1 \\ \mu(n)^2=1}}\frac{\rho_j(n)\varphi(n)}{n^2}\geq\Sigma_1-\Sigma_2.\]
En utilisant l'astuce de Rankin, avec $\alpha:=1/\log(x^{\varepsilon_0/C})$, on a
\begin{align*}
\Sigma_2&\leq x^{-\alpha\varepsilon_0}\sum_{\substack{P^+(n)\leq x^{\varepsilon_0/C} \\ \omega_{y_j}(n)=k_j' \\ (n,W\beta D)=1 \\ \mu(n)^2=1}}\frac{\rho_j(n)\varphi(n)n^\alpha}{n^2}\\
&\leq\e^{-C}\frac{L_j(\alpha)^{k_j'}}{k_j'!}\prod_{\substack{y_j^*<p\leq x^{\varepsilon_0/C} \\ p\nmid W\beta D}}\left(1+\frac{\rho_j(p)(p-1)p^{\alpha}}{p^2}\right)\\
&\ll\e^{-C}\frac{L_j^{k_j'}}{k_j'!}\exp\bigg\{\sum_{\substack{y_j^*<p\leq x \\ p\nmid W\beta D}}\frac{\rho_j(p)}{p}\bigg\}.
\end{align*}
Par ailleurs, on minore $\Sigma_1$ pour $k_j'\geq 1$ par
\begin{align*}
\Sigma_1&\geq \bigg(\frac{L_j(0)^{k_j'}}{k_j'!}+O\Big(\sum_{\substack{w<p\leq y_j^* \\ p\nmid\beta D}}\frac{1}{p^2}\frac{L_j(0)^{k_j'-1}}{(k_j'-1)!}\Big)\bigg)\prod_{\substack{y^*<p\leq x^{\varepsilon_0/C} \\ p\nmid W\beta D}}\left(1+\frac{\rho_j(p)(p-1)}{p^2}\right)\\
&\gg\bigg(\frac{L_j^{k_j'}}{k_j'!}+O\Big(\frac{L_j^{k_j'-1}}{w(k_j'-1)!}\Big)\bigg)\exp\bigg\{\sum_{\substack{y_j^*<p\leq x \\ p\nmid W\beta D}}\frac{\rho_j(p)}{p}\bigg\}C^{-g_j}.
\end{align*}
La même inégalité est valable sans le terme d'erreur dans la parenthèse pour $k_j'=0$. Finalement, puisque $k_j'\leq L_j$, lorsque $w$ et $C$ sont assez grands, on obtient bien le lemme énoncé.
\end{proof}

On peut désormais démontrer~(\ref{vejfhbnqmeflsdj}) pour le choix~(\ref{jhljshgdflh}). Grâce au \og$W$-trick\fg, les conditions $(a_i, a_j)=1$ sont négligeables lorsque $w$ est assez grand. On introduit donc
\[\Sigma_3:=\prod_{1\leq j\leq r}\sum_{\substack{a_j\leq x^{\varepsilon_0} \\ \omega_{y_j}(a_j)=k_j \\ (a_j, W\beta D)=1 \\ \mu(a_j)^2=1}}\frac{\rho_j(a_j)\varphi(a_j)}{a_j^2}\]
et
\begin{multline*}
\Sigma_4:=\sum_{\substack{d_{\mathcal{I}}\leq x^{\varepsilon_0}\ (\mathcal{I}\subseteq\mathcal{P}(\llbracket 1, r\rrbracket), \vert\mathcal{I}\vert\geq 2) \\ \prod_{\mathcal{I}}d_{\mathcal{I}}\neq 1\\ (d_{\mathcal{I}}, W\beta D)=1 \\ \mu(d_{\mathcal{I}})^2=1}}\bigg(\prod_{\mathcal{J}}\frac{\prod_{j\in\mathcal{J}}\rho_j(d_{\mathcal{J}})}{d_{\mathcal{J}}^{\vert\mathcal{J}\vert}}\bigg)\\
\times\prod_{1\leq j\leq r}\sum_{\substack{a_j\leq x^{\varepsilon_0} \\ \omega_{y_j}(a_j)=k_j-\sum_{\mathcal{J}\ni j}\omega_{y_j}(d_{\mathcal{J}}) \\ (a_j, W\beta D)=1 \\ \mu(a_j)^2=1}}\frac{\rho_j(a_j)\varphi(a_j)}{a_j^2},
\end{multline*}
de sorte que 
\[\sum_{\substack{x<n\leq x+y \\ \omega_{y_j}(Q_j(n))=k_j}}1\gg\frac{\e^{-M}y}{(\log x)^r}\left(\Sigma_3-\Sigma_4\right).\]
Dans la définition de $\Sigma_4$, les variables $d_{\mathcal{I}}$ encodent les relations de coprimalité entre les $a_i$ pour $i\in\mathcal{I}$. Dire que les $a_1,\dots,a_r$ ne sont pas tous premiers entre eux est exactement dire $\prod_{\mathcal{I}}d_{\mathcal{I}}\neq 1$. Avec le Lemme~\ref{hejfsndlsdfuh}, on a d'une part,
\[\Sigma_3\gg\prod_{1\leq j\leq r}\frac{L_j^{k_j}}{k_j!}\exp\bigg\{\sum_{\substack{y_j^*<p\leq x \\ p\nmid\beta D}}\frac{\rho_j(p)}{p}\bigg\},\]
et d'autre part, puisque $\rho_j(p)\leq g$ et $k_j\leq L_j$ pour tout $1\leq j\leq r$,
\[\Sigma_4\ll\sum_{\substack{d_{\mathcal{I}}\leq x^{\varepsilon_0}\ (\mathcal{I}\subseteq\mathcal{P}(\llbracket 1, r\rrbracket), \vert\mathcal{I}\vert\geq 2) \\ \prod_{\mathcal{I}}d_{\mathcal{I}}\neq 1\\ (d_{\mathcal{I}}, W\beta D)=1 \\ \mu(d_{\mathcal{I}})^2=1}}\bigg(\prod_{\mathcal{J}}\frac{g^{r\omega(d_{\mathcal{J}})}}{d_{\mathcal{J}}^{\vert\mathcal{J}\vert}}\bigg)\Sigma_3\ll\frac{\Sigma_3}{\sqrt{w}}.\]
Ainsi, lorsque $w$ est suffisamment grand, pour ce choix de $k_1,\dots,k_r$ on obtient avec~(\ref{mj}),
\[\sum_{\substack{x<n\leq x+y \\ \omega_{y_j}(Q_j(n))=k_j}}1\gg\frac{\e^{-M}y}{(\log x)^r}\prod_{1\leq j\leq r}\frac{\varphi_j(\beta D)\e^{-M_j}\log x}{\beta D\sqrt{1+L_j}},\]
et donc le résultat désiré.
\end{proof}

\bibliographystyle{plain-fr}
\bibliography{Bibliographie}

\noindent\bsc{Institut de Math\'ematiques de Jussieu-PRG, Universit\'e Paris Diderot,
Sorbonne Paris Cit\'e, 75013 Paris, France}\\

\noindent\textit{E-mail :} \url{eliegoudout@hotmail.com}

\end{document}